\documentclass{amsart}
\usepackage{amssymb}
\usepackage{amscd}

\newtheorem{Theorem}{Theorem}[section]
\newtheorem{Lemma}[Theorem]{Lemma}
\newtheorem{Corollary}[Theorem]{Corollary}
\newtheorem{Proposition}[Theorem]{Proposition}
\theoremstyle{definition}
\newtheorem{Definition}[Theorem]{Definition}
\newtheorem{Example}[Theorem]{Example}
\theoremstyle{remark}
\newtheorem{Remark}{Remark}

\font\sy=cmsy10

\font\ym=msbm10  
\newcommand{\cC}{{\hbox{\sy C}}}

\newcommand{\cN}{{\hbox{\sy N}}}

\newcommand{\C}{{\text{\ym C}}}

\renewcommand{\a}{\alpha}
\renewcommand{\b}{\beta}

\renewcommand{\d}{\delta}

\newcommand{\e}{\epsilon}

\renewcommand{\i}{\iota}

\renewcommand{\l}{\lambda}

\newcommand{\m}{\mu}

\renewcommand{\o}{\omega}

\newcommand{\p}{\pi}
\renewcommand{\r}{\rho}
\newcommand{\s}{\sigma}
\renewcommand{\t}{\tau}

\newcommand{\x}{\xi}

\newcommand{\y}{\eta}
\newcommand{\z}{\zeta}
\newcommand{\End}{\hbox{\rm End}}
\newcommand{\Hom}{\hbox{\rm Hom}}

\newcommand{\supp}{\hbox{\rm supp}}

\title[Free Product Bimodules]
{C*-Tensor Categories and Free Product Bimodules}
\author[Yamagami Shigeru]{Shigeru Yamagami\\
Department of Mathematical Sciences\\
Ibaraki University\\
Mito, 310--8512, JAPAN}
\begin{document}
\maketitle
\begin{abstract}
A C*-tensor category with simple unit object is realized 
by von Neumann algebra bimodules of finite Jones index 
if and only if it is rigid. 
\end{abstract}
\bigskip
\noindent{\bf Introduction} 

Given a subfactor $N \subset M$ of finite Jones index, 
the associated sequence of higher relative commutants 
$\{ N'\cap M_j \}_{j \geq 0}$ with $M_0 = M$ contains 
rich information on the relative position of the inclusion 
and has been a good source of combinatorial structures behind
subfactors. 

An axiomatization of higher relative commutants is performed by
S.~Popa in the form of so-called standard lattices, 
which is based on a generalization of his preceding result on 
free product construction of irreducible subfactors of 
arbitrarily given indices. 
On the other hand, F.~Radulescu invented a method to construct 
free product subfactors from commuting squares satisfying some 
strong conditions on non-degeneracy and connectedness. 

Under the background of these results on subfactors, we shall 
present here an abstract characterization of C*-tensor categories 
which can be realized by bimodules of finite Jones index. 

To do this, we need to impose some structure of rigidity in 
C*-tensor categories: we assume Frobenius duality together with 
a systematic choice of dual objects in C*-tensor categories. 
The notion of Frobenius duality is formulated in \cite{Y1} 
as an abstraction from the Ocneanu's Frobenius reciprocity 
in bimodules of finite Jones index (see \cite[Chapter 9]{EK}), 
which produces promptly all the combinatorial features in subfactor 
theory such as commuting squares, Markov traces, Perron-Frobenius 
eigenvectors and so on (\cite{Y4}).

Frobenius duality, together with the associated cyclic tensor
products, is utilized in our previous paper \cite{HY} to realize 
amenable C*-tensor categories as bimodules over the AFD 
$\text{II}_1$-factor, where random walks on the fusion algebra 
(the Grothendieck ring of a tensor category) is coupled with 
the structure of tensor categories. 
The constructed bimodules are referred to as random walk bimodules 
in what follows. 

Note here that random walk bimodules are based on non-factor 
von Neumann algebras and far from being irreducible generally. 
(The keypoint in \cite{HY} is that amenability is enough to 
prove the irreducibility.)

Our main result in this paper 
is that the random walk construction is combined 
with the Radulescu's method to produce bimodules over 
amalgamated free product factors so that it gives a fully 
faithful realization of a C*-tensor category with 
Frobenius duality as bimodules of finite Jones index. 

Since a C*-tensor category with Frobenius duality is 
characterized abstractly as a rigid C*-tensor category with 
simple unit object (\cite{Y5}), 
the above result gives the following more satisfactory 
characterization: 
A C*-tensor category with simple unit object is realized as that of 
bimodules of finite Jones index over a factor if and only if 
it is rigid.

\section{Preliminaries}
In this paper, we shall work with C*-tensor categories, which 
we may assume to be strict without loss of generality by coherence
theorem (\cite[Theorem~7.2.1]{M}). 

A {\bf conjugation} in a C*-tensor category $\cC$ is, by definition, 
a conjugate-linear monoidal C*-functor, $X \mapsto \overline X$, 
$\Hom(X,Y) \ni f \mapsto \overline f 
\in \Hom(\overline X, \overline Y)$ with the accompanied conjugate 
multiplicativity $\{ c_{X,Y}: \overline Y \otimes \overline X \to 
\overline{X\otimes Y} \}$ and a natural equivalence 
$\{ d_X: X \to \overline{\overline X} \}$ satisfying 
$\overline{d_X} = d_{\overline X}$. 
The object $\overline X$ is often denoted by $X^*$ 
with the associated contravariant functor defined by 
${}^tf = \overline{f^*} = {\overline f}^*$ for a morphism $f$.

\begin{Definition}
Let $\cC$, $\cC'$ be C*-tensor categories with conjugations.
A monoidal C*-functor $F:\cC \to \cC'$ with multiplicativity 
$m = \{ m_{X,Y}: F(X)\otimes F(Y) \to F(X\otimes Y) \}$ 
is said to be {\bf *-monoidal}
if there is a natural family 
$\{ s_V: F(V^*) \to F(V)^*\}$ of unitaries in $\cC'$ satisfying 
\[
\begin{CD}
F(W^*)\otimes F(V^*) @>{s\otimes s}>> F(W)^*\otimes F(V)^* 
@>{c}>> (F(V)\otimes F(W))^*\\
@V{m}VV @. @AA{{}^tm}A\\
F(W^*\otimes V^*) @>>{F(c)}> F((V\otimes W)^*) 
@>>{s}> F(V\otimes W)^*
\end{CD}\ ,
\]
\[
\begin{CD}
F(V^{**}) @>{s}>> F(V^*)^*\\
@V{F(d^{-1})}VV @AA{{}^ts}A\\
F(V) @>>{d}> F(V)^{**}
\end{CD}\ .
\]
\end{Definition}

\begin{Definition}
Let $F$, $G: \cC \to \cC'$ be monoidal functors.
A natural C*-transformation $\{\varphi_V: F(V) \to G(V)\}$ is called
a {\bf *-monoidal transformation} if 
it is monoidal (multiplicative) and satisfies 
\[
\begin{CD}
F(V^*) @>{s^F}>> F(V)^*\\
@V{\varphi_{V^*}}VV @AA{{}^t\varphi_V}A\\
G(V^*) @>>{s^G}> G(V)^*
\end{CD}\ .
\]

A *-monoidal transformation is called 
a {\bf *-monoidal equivalence} 
if $\varphi_V$ is an isomorphism for each $V$.
\end{Definition}

\begin{Definition}
A *-monoidal functor $F: \cC \to \cC'$ is called 
a {\bf *-monoidal isomorphism} if 
there are a *-monoidal functor 
$G:\cC' \to \cC$ and two *-monoidal equivalences 
\[
\varphi: FG \cong \text{id}_{\cC'},\qquad 
\psi: GF \cong \text{id}_\cC.
\]
Two C*-tensor categories $\cC$, $\cC'$ with conjugations are said to be 
{\bf isomorphic} if there is a 
*-monoidal isomorphism $F: \cC \to \cC'$.
\end{Definition}

A conjugation is said to be {\bf strict} if $\{ c_{X,Y} \}$ and 
$\{ d_X \}$ are identities. 
The coherence theorem can be extended to tensor categories with 
conjugations and we can safely restrict ourselves to strict 
conjugations (cf.~\cite{Y5}). 

\begin{Definition}
A {\bf Frobenius duality} 
in a (strict) C*-tensor category $\cC$ is
a (strict) conjugation together with a family of self-conjugate morphisms 
$\{ \e_X = \overline{\e_X}: X\otimes X^* \to I\}_{X \in  \text{Object}}$ 
satisfying the following conditions 
($I$ being the unit object in $\cC$). 
\begin{enumerate}
\item
(Multiplicativity)
\[
\begin{CD}
X\otimes Y\otimes Y^*\otimes X^* @>{\epsilon_{X\otimes Y}}>> I\\
@V{1\otimes \epsilon_Y\otimes 1}VV @|\\
X\otimes X^* @>>{\epsilon_X}> I
\end{CD}\ .
\]
\item
(Naturality)
For a morphism $f: X \to Y$ in $\cC$, 
\[
\begin{CD}
X\otimes Y^* @>{f\otimes 1}>> Y\otimes Y^*\\
@V{1\otimes {}^tf}VV @VV{\epsilon_Y}V\\
X\otimes X^* @>>{\epsilon_X}> I
\end{CD}\ .
\]
\item
(Faithfulness)
The map 
\[
Hom(X,Y) \ni f \mapsto \epsilon_Y\circ(f\otimes 1) \in 
Hom(X\otimes Y^*,I)
\]
is injective for $X$, $Y \in \text{Object}(\cC)$. 
\item 
(Neutrality) 
For a morphism $f \in \End(X)$, we have 
\[
\e_X(f\otimes 1)\e_X^* = \e_{X^*}(1\otimes f)\e_{X^*}^*. 
\]
\end{enumerate}
\end{Definition}

%\begin{figure}[htbp]
%\includegraphics[width=12cm,clip]{fig1.eps}
%%\includegraphics{fig1.eps}
%\caption[]{}
%\end{figure}

\begin{Example}
Let $N$ be a $\text{II}_1$-factor and consider the tensor category 
$\cC$ of $N$-$N$ bimodules of finite index. 
Then, together with 
the obvious operation of conjugation, $\cC$ admits a canonical 
Frobenius duality $\{ \e_X \}$ defined by 
\[
\e_X: X\otimes_NX^* = L^2(\End(X_N)) 
\ni x(\t\circ E)^{1/2} \mapsto [X]^{1/4} E(x) \t^{1/2}
\]
for an irreducible $X$ ($\t$ is the normalized trace on $N$, 
$E: \End(X_N) \to N$ is the conditional expectation and 
$[X]$ is the Jones index for the inclusion $N \subset \End(X_N)$) 
and then by 
\[
\e_X = \sum_{i=1}^n \e_{X_i} (T_i\otimes \overline T_i) 
\]
for a general $X$, where $X \cong \bigoplus_{i=1}^n X_i$ 
with 
$\{ T_i: X \to X\i \}$ an orthogonal family of coisometries. 
\end{Example}

\begin{Definition}
An object $X$ in a C*-tensor category $\cC$ is said to be 
{\bf rigid} if we can find an object $Y$ and morphisms 
$\e: X \otimes Y \to I$, $\d: I \to Y\otimes X$ such that 
\[
\begin{CD}
X @>{1\otimes \d}>> X\otimes Y\otimes X @>{\e\otimes 1}>> X
\end{CD}, 
\qquad 
\begin{CD}
Y @>{\d\otimes 1}>> Y\otimes X\otimes Y @>{1\otimes \e}>> Y
\end{CD}
\]
are identities. 

An object $Y$ in the definition of rigidity turns out to be 
unique up to isomorphism and is referred to as a {\bf dual object} of 
$X$ (see \cite{K} for example).

A C*-tensor tensor category $\cC$ is {\bf rigid} if any object in 
$\cC$ is rigid. 
\end{Definition}

It is easy to see that a C*-tensor category is rigid if it admits 
a Frobenius duality. Conversely, we have the following by \cite{Y5}. 

\begin{Theorem}
In a rigid C*-tensor category $\cC$ with simple unit object, 
there exists a Frobenius duality and any Frobenoius duality in 
$\cC$ is unique up to unitary isomorphisms. 
\end{Theorem}

The next result is due to Longo and Roberts (\cite[Lemma~3.2]{LR}). 
We here present an independant structural proof. 

\begin{Proposition}[Longo-Roberts]
Let $\cC$ be a C*-tensor category with simple unit object. 
Then for any rigid object $X$ in $\cC$, 
$\End(X)$ is finite-dimensional. 

In particular, the tensor category $\cC$ is semisimple, i.e., 
any object $X$ in $\cC$ is isomorphic to a direct sum of 
finitely many simple objects (by adding subobjects to $\cC$ 
if necessary). 
\end{Proposition}

\begin{proof}
Let $\e: X^*\otimes X \to I$, $\delta: I \to X\otimes X^*$ 
be a rigidity pair and $F: \End(X) \to \Hom(I, X\otimes X^*)$ 
be the associated Frobenius transform: 
$F(f) = (f\otimes 1)\delta$ and 
$F^{-1}(g) = (1\otimes \e)(g\otimes 1)$. 
Note here that $\Hom(I,X\otimes X^*)$ is a Hilbert space 
because of the simplicity of the unit object $I$.

From the inequalities
\begin{align*}
\| F(f) \|^2 &= \delta^*(f^*f\otimes 1)\delta 
\leq \| f\|^2 \| \delta \|^2,\\
\| F^{-1}(g) \|^2 &= (1\otimes \e)(gg^*\otimes 1)(1\otimes \e^*) 
\leq \| g\|^2 \| \e\|^2,
\end{align*}
the C*-algebra $\End(X)$ is continuously 
isomorphic to the Hilbert space $\Hom(I,X\otimes X^*)$ 
with the bounded inverse, whence 
$\End(X)$ is reflexive as a Banach space, proving 
$\dim(\End(X)) < +\infty$. 
\end{proof}

\begin{Proposition}
Let $N$ be a factor. An $N$-$N$ bimodule $X$ has finite Jones index 
if and only if $X$ is rigid in the tensor category $\cC$ 
of $N$-$N$ bimodules.
\end{Proposition}

\begin{proof}
If $X$ has the finite Jones index, it is rigid as a consequence of 
the existence of Frobenius duality. 

Conversely, let $X$ be a rigid object in $\cC$. Then by Longo-Roberts' 
finite-dimensionality, we may assume that $X$ is simple. 
In that case, the finiteness of Jones index follows from the
non-triviality of 
$\Hom({}_NL^2(N)_N, {}_NX\otimes_NX^*_N)$ and 
$\Hom({}_NL^2(N)_N, {}_NX^*\otimes_NX_N)$ as pointed out in 
\cite[Lemma~10]{KY}.
\end{proof}

In what follows, we shall work with a rigid C*-tensor category $\cC$ 
with simple unit object (Frobenius duality being implicitly assumed 
by Theorem~1.7). 
Let $S$ be the set of equivalence classes of simple objects in $\cC$. 
By the semisimplicity of $\cC$, the free vector space $\C[S]$ 
over the set $S$ admits the algebra structure as a Grothendieck ring. 
Moreover, the existence of dual objects in $\cC$ allows us to define 
the *-operation in $\C[S]$ by $[X]^* = [X^*]$. 
The *-algebra $\C[S]$ is, by definition, the {\bf fusion algebra} 
of $\cC$. 

Now we review the construction and some of its basic properties 
introduced in \cite{HY}.
By choosing a representative set of simple objects in $\cC$, 
we regard $S$ as a set of simple objects in $\cC$. 

Let $\m$ be a probability measure on $S$ which supports 
the whole set $S$. To introduce the random walk construction of 
bimodules, we choose an orthogonal family $\{ e_s \}_{s \in S}$ 
of projections in the AFD $\text{II}_1$-factor $R$ such that 
\[
\o(e_s) = \frac{\m(s)}{d(s)}, 
\]
where $\o$ denotes the normalized trace of $R$. 

For a finite sequence $x = (x_n,\dots,x_1) \in S^n$, 
we define the projection $e_x \in R^{\otimes n}$ by 
$e_x = e_{x_n}\otimes \dots \otimes e_{x_1}$. 

We denote the space $e_xR^{\otimes n}e_y$ simply by ${}_xR_y$ 
and keep to use $\o$ to denote 
the normalized trace on $R^{\otimes n}$. 

In the following, 
the conventional symbol $\otimes$ for tensor products 
is often omited to simplify notations. 

Now, given an object $X$ in $\cC$, we define an increasing 
sequence of finite von Neumann algebras 
$\{ A_n(X) \}_{n \geq 0}$ and an increasing sequence of 
Hilbert spaces $\{ X_n \}_{n \geq 0}$ by 

\begin{align*}
A_n(X) &= \bigoplus_{x, y \in S^n} 
\begin{bmatrix}
x_n\dots x_1X\\
y_n\dots y_1X
\end{bmatrix}
\otimes {}_xR_y, \\
X_n &= \bigoplus_{x, y \in S^n} 
\begin{bmatrix}
x_n\dots x_1X\\ y_n\dots y_1
\end{bmatrix}
\otimes L^2({}_xR_y),
\end{align*}
with the imbeddings 
$A_n(X) \to A_{n+1}(X)$ and $X_n \to X_{n+1}$ defined by 
\begin{align*}
\s\otimes a &\mapsto \sum_{s \in S} 
(1_s\otimes \s)\otimes (e_s\otimes a),\\
\x\otimes a\o^{1/2} &\mapsto \sum_{s \in S}
(1_s\otimes\x)\otimes (e_s\otimes a)\o^{1/2}.
\end{align*}

If we introduce a faithful tracial functional $\t_X^n$ on $A_n(X)$ by 
\[
\t_X^n(\s\otimes a) = \d_{x,y} \langle \s \rangle \o(a)
\quad\text{with $\langle \s \rangle 1_I = \e_{xX}(\s\otimes 1)\e_{xX}^*$,}
\]
then it turns out that the family $\{ \t_X^n \}$ is compatible with 
the inclusion $A_n(X) \to A_{n+1}(X)$, whence it defines the trace 
$\t_X$ on the inductive limit $\cup_{n \geq 0} A_n(X)$. 

The standard Hilbert space $L^2(A_n(X))$ of $A_n(X)$ is naturally
identified with $(XX^*)_n$ by 
\[
(\s\otimes a)\t_X^{1/2} \leftrightarrow \widetilde\s \otimes
a\o^{1/2}, 
\]
where 
$\displaystyle \widetilde \s \in 
\begin{bmatrix}
x_n\dots x_1XX^*\\ y_n\dots y_1
\end{bmatrix}$ 
is the Frobenius transform of $\s$. 

If we write $A_n$ to stand 
for $A_n(I)$ with the trace $\t_I$ denoted by $\t$, 
the Hilbert space $X_n$ is an $A_n$-$A_n$ bimodule by 
\[
(\s\otimes a)(\x\otimes x\o^{1/2})(\s'\otimes a') 
= (\s\otimes 1_X)\x \s' \otimes axa'\o^{1/2}, 
\]
which is referred to as a {\bf random walk bimodule}.

For later use, we present here the following fact, 
which is immediate from definitions 
(the details can be found before Proposition~3.8 in \cite{HY}).

\begin{Lemma}
Every algebraic vector in $X_n$ is $\t_n$-bounded and 
the associated $A_n$-valued inner product is given by the formula
\[
{}_{A_n}[\x\otimes a\o^{1/2}, \y\otimes b\o^{1/2}] 
= \langle \x\y^*\rangle_X \otimes ab^* \in A_n, 
\]
where 
$\displaystyle \x \in 
\begin{bmatrix}
xX\\ x'
\end{bmatrix}$, 
$\displaystyle \y \in 
\begin{bmatrix}
yX\\ y'
\end{bmatrix}$, 
$a \in {}_xR_{x'}$, $b \in {}_yR_{y'}$ and 
\[
\langle \x\y^* \rangle_X 
= (1\otimes \e_X)(\x\y^*\otimes 1_{X^*})(1\otimes \e_X^*) 
\in 
\begin{bmatrix}
x\\ y
\end{bmatrix}
\]
denotes the partial trace of $\x\y^*$. 
\end{Lemma}

%\begin{figure}[htbp]
%\includegraphics{fig2.eps}
%\vspace{0.5cm}
%\caption[]{}
%\end{figure}

\section{Free Product Bimodules}

Among finite von Neumann algebras constructed by random walks, 
we here concentrate on the lowest inclusion $A_1(X) \subset A_2(X)$ 
(higher inclusions can be used as well). 
To simplify the notation, we denote these as $A(X)$ and $B(X)$ with 
the convention $A = A(I)$ and $B = B(I)$. 

Recall that all these are isomorphic to a direct sum of countably many 
AFD $\text{II}_1$-factors. 

We also use the notation ${}_AX_A$ and ${}_BX_B$ to stand for 
random walk bimodules $X_1$ and $X_2$ resepectively. 
(The notation in fact indicates the fact that ${}_AX_A$ is an $A$-$A$
bimodule.)

\begin{Lemma}\label{C}
The inclusion $A \subset B$ is connected, i.e., 
\[
Z(A) \cap Z(B) = \C 1.
\]
\end{Lemma}

\begin{proof}
Recall that 
\[
A = \bigoplus_{s,t \in S} 
\begin{bmatrix}
s\\ t
\end{bmatrix}
\otimes {}_sR_t 
= \bigoplus_{s \in S} 
\begin{bmatrix}
s\\ s
\end{bmatrix}
\otimes {}_sR_s 
\cong \bigoplus_{s \in S} {}_sR_s
\]
is imbedded into 
\[
B = \bigoplus_{s, t \in S^2} 
\begin{bmatrix}
s\\ t
\end{bmatrix}
\otimes {}_sR_t
\]
by 
\[
1_s\otimes x \mapsto \sum_{t \in S} (1_t\otimes 1_s)\otimes 
(e_t\otimes x).
\]
Thus a central element $\sum_s f(s) 1_s\otimes e_s$ in $A$ is 
mapped into 
\[
\sum_{s,t} f(s) 1_{st}\otimes e_{st} = 
\sum_u \sum_{st} \sum_{\x: u \to st} f(s)
\frac{d(u)}{(\x|\x)} \x\x^*\otimes e_{st}. 
\]
If we identify this with a central element 
\[
\sum_u g(u) \sum_{s,t} \sum_{\x: u \to st} 
\frac{d(u)}{(\x|\x)} \x\x^*\otimes e_{st}
\]
in $B$, then we have 
\[
f(s) \sum_{\x: u \to st} \frac{d(u)}{(\x|\x)} \x\x^* 
= g(u) \sum_{\x: u \to st} \frac{d(u)}{(\x|\x)} \x\x^*
\]
for any $s$, $t$ and $u \in S$ ($e_{st} \not= 0$ as $\m$ supports 
the whole $S$), i.e., 
\[
f(s) = g(u) \quad\text{for any $s$, $u \in S$}
\]
because 
$\displaystyle 
\begin{bmatrix}
st\\ u
\end{bmatrix}
\not= \{ 0\}$ for some $t \in S$. 
\end{proof}

\begin{Lemma}
For $s$, $t$, $x$, $y$ and $z \in S$, 
\[
{}_{xy}R_{sz}(e_s\otimes {}_zR_t) = {}_{xy}R_{st}.
\]
\end{Lemma}

\begin{proof}
Consider the case $\m(z)/d(z) \leq \m(t)/\d(t)$. 
By the assumption $\m(z) > 0$ for any $z \in S$, 
we can find a finite family of partial isometries 
$\{ u_i \}_{0 \leq i \leq n}$ in $R$ 
($n$ can be chosen as the positive integer satisfying 
$n \leq \m(t)d(z)/d(t)\m(z) < n+1$) such that 
\[
e_t = \sum_i u_i^*e_zu_i.
\]
Then $u_ie_t \in {}_zR_t$ and we have 
\begin{align*}
{}_{xy}R_{sz}(e_s\otimes {}_zR_t) 
&\supset \sum_i {}_{xy}R_{sz} (e_s\otimes u_ie_t)\\
&\supset \sum_i e_{xy} (R\otimes R) (e_s\otimes u_i^*e_zu_ie_t)\\
&= {}_{xy}R_{st}.
\end{align*}
The reverse inclusion is trivial.
\end{proof}

\begin{Lemma}\label{D}
Both of ${}_AX_AB$ and $B{}_AX_A$ are dense in ${}_BX_B$.
\end{Lemma}

\begin{proof}
We will only show the density of $B{}_AX_A$. 
Let $\x \in {}_BX_B$ be orthogonal to $B{}_AX_A$ and express it as 
\[
\x = \bigoplus_{s,t \in S^2} \sum_i \x_i(s,t)\otimes \a_i(s,t)
\]
with $\{\a_i(s,t)\}_{i \geq 1}$ an orthonormal basis in $L^2({}_sR_t)$ 
and $\x_i(s,t) \in \Hom(t,sX)$.

Given any $u \in S$, choose $\s \in \Hom(t_2u,s)$, 
$\y \in \Hom(t_1,uX)$, $a \in {}_sR_{t_2u}$ and $b \in {}_uR_{t_1}$ 
arbitrarily. Then $\x$ is orthogonal to the vector 
\[
(\s\otimes 1_X)(1_{t_2}\otimes\y) \otimes a(e_{t_2}\otimes b)\o^{1/2}
\]
by assumption, i.e., 
\[
\sum_i (\x_i(s,t)|(\s\otimes 1_X)(1_{t_2}\otimes\y))\, 
(\a_i(s,t)|a(e_{t_2}\otimes b)\o^{1/2}) = 0. 
\]
Since the set $\{ a(e_{t_2}\otimes b)\o^{1/2} \}$ is total in 
$L^2({}_sR_t)$ by Lemma~2.2, we can simplify the orthogonality 
condition to 
\[
(\x_i(s,t)|(\s\otimes 1_X)(1_{t_2}\otimes\y)) = 0 
\]
for any $i \geq 1$. 

On the other hand, Frobenius isomorphisms 
\[
\bigoplus_{u \in S} 
\begin{bmatrix}
s\\ t_2u
\end{bmatrix}
\otimes
\begin{bmatrix}
uX\\ t_1
\end{bmatrix}
\cong 
\bigoplus_{u \in S} 
\begin{bmatrix}
t_2^*s\\ u
\end{bmatrix}
\otimes 
\begin{bmatrix}
u\\ t_1X^*
\end{bmatrix}
\cong 
\begin{bmatrix}
t_2^*s\\ t_1X^*
\end{bmatrix}
\cong 
\begin{bmatrix}
sX\\ t
\end{bmatrix}
\]
show that $\{ (\s\otimes 1_X)(1_{t_2}\otimes \y) \}$ is total 
in $\Hom(t,sX)$. 

Consequently we have $\x_i(s,t) = 0$ for any $i \geq 1$ and any 
$s$, $t \in S^2$. 
\end{proof}

\begin{Lemma}\label{Shift}
Given an object $X$ of $\cC$, the following formulas define unitary maps
$L^2(B)\otimes_A X_A \to {}_BX_B$ and 
${}_AX\otimes_AL^2(B) \to {}_BX_B$:
\begin{gather*}
(\s\otimes b\o^{1/2}) \otimes_A (\x\otimes x\o^{1/2}) 
\mapsto 
\sum_{u \in S} (\s\otimes 1_X)(1_u\otimes \x) \otimes 
b(e_u\otimes x)\o^{1/2}, 
\\
(\x\otimes x\o^{1/2})\otimes_A (\s\otimes b\o^{1/2}) 
\mapsto 
\sum_{u \in S} (1_u\otimes \x)\s \otimes (e_u\otimes x)b \o^{1/2},
\end{gather*}
where 
\[
\s \in 
\begin{bmatrix}
s_2s_1\\ t_2t_1
\end{bmatrix},\  
b \in {}_{s_2s_1}R_{t_2t_1}\quad\text{and}
\quad 
\x \in 
\begin{bmatrix}
sX\\ t
\end{bmatrix}, \ 
x \in {}_sR_t
\]
with $s$, $t$, $s_i$, $t_i \in S$.
\end{Lemma}

\begin{proof}
By the formula of $A$-valued inner product in ${}_AX$ 
(Lemma~1.10), 
\begin{align*}
\| (\s\otimes b\o^{1/2})\otimes_A (\x\otimes x\o^{1/2}) \|^2 
&= \sum_u (\s\otimes b\o^{1/2}| 
(\s\otimes b\o^{1/2}) (1_u\otimes \langle \x\x^*\rangle_X 
\otimes (e_u\otimes xx^*)))\\
&= \sum_u \langle (\s^*\s\otimes 1_X)(1_u\otimes \x\x^*) \rangle 
\o(b^*b(e_u\otimes xx^*))\\
&= \| \sum_u (\s\otimes 1_X)(1_u\otimes \x) \otimes 
b(e_u\otimes x)\o^{1/2} \|^2,
\end{align*}
i.e., the correspondance 
\[
(\s\otimes b\o^{1/2}) \otimes_A (\x\otimes x\o^{1/2}) 
\mapsto 
\sum_{u \in S} (\s\otimes 1_X)(1_u\otimes \x) \otimes 
b(e_u\otimes x)\o^{1/2} 
\]
gives a well-defined isometry $L^2(B)\otimes_AX_A \to {}_BX_B$, 
which is in fact a unitary by Lemma~\ref{D}.
\end{proof}

The unitary maps in the above lemma are clearly functorial in 
the variable $X$ and, composing these unitary maps, we obtain 
a natural family of unitary maps 
\[
{}_AL^2(B)\otimes_A X_A \to {}_AX\otimes_A L^2(B)_A, 
\]
which intertwines the $A$-$A$ actions and is referred to as 
{\bf shift isomorphisms}. 

\begin{Corollary}
The shift isomorphism transfers the subspace $L^2(B)^\circ\otimes_AX_A$ 
onto the subspace ${}_AX\otimes_AL^2(B)^\circ$. Here 
$L^2(B)^\circ$ denotes the orthognal complement of $L^2(A)$ in 
$L^2(B)$ ($L^2(A)$ being imbedded into $L^2(B)$ by the
trace-preserving conditional expectation).  
\end{Corollary}

\begin{proof}
From the definition of unitaries, it is immediate to check 
that both of $\t_B^{1/2}\otimes_{\t_A^{1/2}} \x$ and 
$\x\otimes_{\t_A^{1/2}} \t_B^{1/2}$ are mapped to $\x \in {}_AX_A$. 
%\[
%\t_B^{1/2} = \sum_{s \in S^2} 1_s\otimes e_s\o^{1/2}
%\]
Thus the unitaries turn out to be 
the obvious identification 
\[
L^2(A)\otimes_A X_A = {}_AX_A = {}_AX\otimes_AL^2(A)
\]
on the subspaces $L^2(A)\otimes_AX_A$ and ${}_AX\otimes_AL^2(A)$. 
Since $L^2(B)^\circ\otimes_AX_A$ and ${}_AX\otimes_AL^2(B)^\circ$ 
are orthogonal complements of these, the assertion holds. 
\end{proof}

We are now ready to introduce a free product construction of 
bimodules. 
Given a continuous von Neumann algebra $Q$ with a faithful 
normalized trace $\t_Q$, we set 
\[
N = N(Q) = (Q\otimes A)*_AB, 
\]
which is a $\text{II}_1$-factor by the connectedness of 
the inclusion $A \subset B$ (see \cite[\S 1]{R} or Lemma~3.8 below).

Given an object $X$ in $\cC$, we shall construct an $N$-$N$ bimodule 
in the following way. 
To define the base Hilbert space of a free-product bimodule, 
we start with the Hilbert space ${}_AX_A$
%\[
%L^2(Q)\otimes {}_AX_A = {}_AX_A\otimes L^2(Q)
%\]
%(the equality means $Q$ commutes with ${}_AX_A$ completely) 
and take $A$-tensor products with $A$-$A$ bimodules
\[
L^2(Q)^\circ\otimes L^2(A) = L^2(A)\otimes L^2(Q)^\circ 
\quad\text{and}\quad 
L^2(B)^\circ = L^2(B) \ominus L^2(A)
\]
alternately so that 
\begin{enumerate}
\item
tensor products of the following two types are not allowed to appear, 
\[
\bigl( L^2(Q)^\circ\otimes L^2(A) \bigr) \otimes_AX_A\otimes_A 
\bigl( L^2(Q)^\circ \otimes L^2(A) \bigr), 
\quad 
L^2(B)^\circ \otimes_A X_A \otimes_A L^2(B)^\circ,
\]
\item 
$L^2(Q)^\circ$ commutes with $L^2(A)A$ and ${}_AX_A$, 
\item
if the tensor component of the form $L^2(B)^\circ\otimes_AX$
appears, it is identified with 
$X\otimes_AL^2(B)^\circ$ by the shift isomorphism. 
\end{enumerate}
Note that tensor components can be lined up sequentially 
for the identification in (iii) and hence there arise no coherence 
problems here.

The direct sum of resulting Hilbert spaces is denoted by $F(X)$, 
which we shall make into an $N$-$N$ bimodule. 

By the above requirements for identification, 
the position of ${}_AX_A$ can be freely moved inside summands. 
For example, 
if we move the component ${}_AX_A$ to the right end in each summand, 
we have the following expression for $F(X)$: 
\begin{align*}
F(X) &= 
{}_AX_A \oplus 
\bigl( L^2(Q)^\circ \otimes L^2(A) \bigr) \otimes_A X_A 
\oplus L^2(B)^\circ \otimes_A X_A\\
&\qquad \oplus L^2(B)^\circ \otimes_A 
\bigl( L^2(Q)\otimes L^2(A) \bigr) \otimes_A {}_AX_A\\
&\qquad \oplus 
(L^2(Q)^\circ \otimes L^2(A)) \otimes_A L^2(B)^\circ \otimes_A
{}_AX_A\\
&\qquad \oplus 
\big( L^2(Q)^\circ \otimes L^2(A) \bigr) 
\otimes_A L^2(B)^\circ \otimes_A 
\bigl( L^2(Q)^\circ \otimes L^2(A) \bigr) 
\otimes_A X_A\\
&\qquad \oplus 
L^2(B)^\circ \otimes_A (L^2(Q)^\circ \otimes L^2(A)) \otimes_A 
L^2(B)^\circ \otimes_A X_A\\
&\qquad \oplus \dots,
\end{align*}
which is further reduced to 
\begin{align*}
F(X) &= {}_AX_A + L^2(Q)^\circ \otimes {}_AX_A 
+ L^2(B)^\circ \otimes_A X_A\\
&\qquad + L^2(B)^\circ \otimes_A L^2(Q)^\circ \otimes {}_AX_A\\
&\qquad + L^2(Q)^\circ \otimes L^2(B)^\circ \otimes_A X_A\\
&\qquad + L^2(Q)^\circ \otimes L^2(B)^\circ \otimes_A 
L^2(Q)^\circ \otimes {}_AX_A\\
&\qquad + L^2(B)^\circ \otimes_A L^2(Q)^\circ \otimes_A 
L^2(B)^\circ \otimes {}_AX_A\\
&\qquad + \dots
\end{align*}
if we use the obvious identification 
$L^2(A)\otimes_AV = V = V\otimes_AL^2(A)$ for an $A$-$A$ bimodule
$V$. 

Let $\cN$ be the dense *-subalgebra of $N = (Q\otimes A)*_AB$ 
algebraically generated by $Q \cup B$. 
If we rearrange $F(X)$ so that the first factors in summands are 
coupled to $L^2(Q) = \C \oplus L^2(Q)^\circ$, then we obtain the expression 
\begin{align*}
F(X) &= L^2(Q)\otimes {}_AX_A + L^2(Q)\otimes
L^2(B)^\circ\otimes_AX_A\\ 
&\qquad+ L^2(Q)\otimes L^2(B)^\circ\otimes_AL^2(Q)^\circ\otimes_AX_A 
+ \dots, 
\end{align*}
on which $Q$ acts from left by multiplication. 
Let us denote this by $\l_Q$. 
On the other hand, if we consider the rearrangement to form 
$L^2(B) = L^2(A) \oplus L^2(B)^\circ$ 
at the left ends of summands, we get the expression like 
\begin{align*}
F(X) &= L^2(B)\otimes_AX_A 
+ L^2(B)\otimes_A L^2(Q)^\circ\otimes_AX_A\\ 
&\qquad+ L^2(B)\otimes_A L^2(Q)^\circ \otimes_A L^2(B)^\circ \otimes_A X_A 
+ \dots
\end{align*}
with the left representation of $B$ denoted by $\l_B$. 
It is then immediate to check the commutativity of $\l_Q(Q)$ 
and $\l_B(A)$, whence we have a *-representation $\l$ of $\cN$ on
$F(X)$ so that $\l|_Q = \l_Q$ and $\l|_B = \l_B$.  

Similarly, we can define an antirepresentation $\rho$ of $\cN$ on
$F(X)$ by right multiplication after rearrangements of summands 
at the right ends. 
It then turns out that $\l$ and $\rho$ commute by straightforward 
computations and the Hilbert space $F(X)$ becomes an $\cN$-$\cN$
module. 

Note that by the way of constructions, the $\cN$-$\cN$ action 
is described according to the free product prescription (cf.~\cite{VDN}). 

\begin{Lemma}
Let $K$ be an object in $\cC$ and $V$ be an $A(K)$-$A(K)$ module. 
Then we can define $A(K)$-$A(K)$ linear unitary maps 
$L^2(B)^\circ\otimes_AV \to L^2(B(K))^\circ \otimes_{A(K)}V$ 
and $V\otimes_A L^2(B)^\circ \to V\otimes_{A(K)} L^2(B(K))^\circ$ 
so that 
\begin{gather*}
(\s\otimes b)\t_B^{1/2} \otimes_{\t_A^{-1/2}} v 
\mapsto (\s\otimes b)\t_K^{1/2} \otimes_{\t_K^{-1/2}} v, \\
v\otimes_{\t_A^{-1/2}} \t_B^{1/2}(\s\otimes b) 
\mapsto v\otimes_{\t_K^{-1/2}} \t_K^{1/2} (\s\otimes b),
\end{gather*}
where $\s \in \Hom(t,s)$ and $b \in {}_sR_t$ with $s$, $t \in S^2$. 
\end{Lemma}

\begin{proof}
By Lemma~\ref{Shift}, 
\[
V\otimes_{A(K)} L^2(B(K))^\circ \cong 
V\otimes_{A(K)} L^2(A(K))\otimes_A L^2(B)^\circ 
\cong V\otimes_A L^2(B)^\circ, 
\]
where an element $v\otimes_{\t^{-1/2}} (\s\otimes b)\t^{1/2}$ 
in the right end is transfered as 
\begin{align*}
v\otimes_{\t_A^{-1/2}} (\s\otimes b)\t_B^{1/2} 
&\leftrightarrow v \otimes_{\t_K^{-1/2}} \t_K^{1/2} 
\otimes_{\t_A^{-1/2}} (\s\otimes b)\t_B^{1/2}\\
&\leftrightarrow v \otimes_{\t_K^{-1/2}} \t_K^{1/2}(\s\otimes b). 
\end{align*}
\end{proof}

\begin{Remark}
A vector $(\s\otimes b)\t_K^{1/2} = \t_K^{1/2}(\s\otimes b)$ in 
$L^2(B(K))$ corresponds to 
the vector $(\s\otimes 1_K)\otimes b\o^{1/2}$ in ${}_BKK^*_B$.
\end{Remark}

\begin{Lemma}
The obvious imbedding $\cN \to (Q\otimes A(K))*_{A(K)}B(K)$ 
is extended to an injective normal *-homomorphism of $N$. 
\end{Lemma}

\begin{proof}
Let $K$ be an object in $\cC$. Then we have the following 
commuting squares 
\[
\begin{CD}
Q\otimes A(K) @>>> A(K)\\
@VVV @VVV\\
Q\otimes A @>>> A
\end{CD},
\qquad
\begin{CD}
A(K) @<<< B(K)\\
@VVV @VVV\\
A @<<< B
\end{CD}
\]
and then we can apply the imbedding theorem of amalgamated free products 
(see Lemma~6.1 in \cite{D} for example) 
to find that $N = (Q\otimes A)*_AB$ is 
a von Neumann subalgebra of the amalgamated free product 
$(Q\otimes A(K))*_{A(K)}B(K)$. 
\end{proof}

\begin{Lemma}
The Hilbert space $F(K\otimes K^*)$ is naturally isometrically isomorphic to 
$L^2\bigl( (Q\otimes A(K))*_{A(K)} B(K) \bigr)$ with the $\cN$-$\cN$
action on $F(K\otimes K^*)$ identified with the left and right action of 
$\cN$ by multiplication when $N$ is imbedded into 
$(Q\otimes A(K))*_{A(K)}B(K)$ by the previous lemma.
\end{Lemma}

\begin{proof}
Let $M = (Q\otimes A(K))*_{A(K)}B(K)$. Then we have 
\[
L^2(M) = L^2(A(K)) \oplus L^2(Q)^\circ\otimes L^2(A(K)) 
\oplus L^2(B(K))^\circ \oplus \dots, 
\]
where the summation is taken over alternate $A(K)$-tensor products of 
$L^2(Q)^\circ\otimes L^2(A(K))$ and $L^2(B(K))^\circ$. 

If we apply the isomorphism in Lemma~\ref{Shift} from outside 
in each summand 
(note that $L^2(Q)^\circ$ commutes with $L^2(A(K))$ and does not touch
on the $A(K)$-action), we end up with a realization of a summand 
in $F(KK^*)$ with the tensor component ${}_AKK^*_A$ at a 
prescribed position: 
two realizations with ${}_AKK^*_A$ at adjacent positions correspond to 
the choices of isomorphisms 
$L^2(B(K))^\circ \cong L^2(B)^\circ\otimes_A L^2(A(K))$ or 
$L^2(B(K))^\circ \cong L^2(A(K))\otimes_A L^2(B)^\circ$
at the last stage of reductions. 

There are two ways of ambiguity other than the ones related by 
shift isomorphisms as indicated by the following diagram 

\[
\begin{CD}
L^2(B(K))^\circ \otimes_{A(K)} L^2(B(K))^\circ 
@>>> L^2(B)^\circ\otimes_A L^2(A(K))\otimes_{A(K)} L^2(B(K))^\circ\\
@VVV @VVV\\
L^2(B(K))^\circ \otimes_{A(K)} L^2(A(K))\otimes_A L^2(B)^\circ 
@.  
L^2(B)^\circ\otimes_A L^2(B(K))^\circ\\
@VVV @VVV\\
L^2(B(K))^\circ\otimes_A L^2(B)^\circ 
@>>> L^2(B)^\circ \otimes_A L^2(A(K))\otimes_A L^2(B)^\circ.
\end{CD}
\]
This is however nothing but two ways of reductions 
in the right hand side of 
\[
L^2(B(K))^\circ \otimes_{A(K)} L^2(B(K))^\circ 
\cong 
\bigl( L^2(B)^\circ\otimes_A L^2(A(K)) \bigr) 
\otimes_{A(K)} 
\bigl( L^2(A(K))\otimes_A L^2(B)^\circ \bigr)
\]
and the commutativity of the diagram is eventually reduced to 
the associativity of multiplication in the algebra $A(K)$. 

%The structure of isomorphisms is visualized by Fig.~\ref{wheel}, 
%where the central big circle denotes a summand in $L^2(M)$ and 
%peripheral dots refer to various realizations 
%(according to the position of ${}_AKK^*_A$) of the corresponding 
%block in $F(KK^*)$, with rays denoting the well-defined isomophism 
%described above and peripheral edges representing 
%(amplified) shift isomorphisms. 

%Note that the commutativity of each triangle is nothing but the
%definition of shift isomorphisms in view of the discussions so far. 

It is now straightforward to check that the $\cN$-$\cN$ action on 
$F(KK^*)$ corresponds to the left and right multiplications of 
$\cN$ on $L^2(M)$. 
\end{proof}

%\begin{figure}[htbp]
%\includegraphics{fig3.eps}
%\caption[]{}
%\label{wheel}
%\end{figure}

\begin{Proposition}
The $\cN$-$\cN$ action on $F(X)$ is extended to the $N$-$N$ action 
by weak continuity: the Hilbert space $F(X)$ is furnished with 
the structure of an $N$-$N$ bimodule. 
\end{Proposition}

\begin{proof}
Let $K$ be such that $X \subset KK^*$
($K = I \oplus X$ for example).
Then $F(X)$ is realized as an $\cN$-$\cN$ invariant closed subspace of 
$F(KK^*)$. 
Since the $\cN$-$\cN$ action is continuously extended to the 
$N$-$N$ action on $F(KK^*)$ by previous lemmas, 
the same holds on $F(X)$. 
\end{proof}

It is also clear that a morphism $f: X \to Y$ in $\cC$ induces 
an $N$-$N$ intertwiner $F(f): F(X) \to F(Y)$ so that 
\[
F(f)(\x\otimes x\o^{1/2}) = (1_s\otimes f)\x \otimes x\o^{1/2}
\]
for $\x \in \Hom(t,sX)$ and $x \in {}_sR_t$ with $s$, $t \in S$. 

The correspondance $X \mapsto F(X)$, together with 
$f \mapsto F(f)$, gives a C*-functor from $\cC$ into the 
C*-tensor category of $N$-$N$ bimodules, which is referred to as 
a {\bf free product functor} in what follows. 

\begin{Lemma}
The free product functor is faithful, i.e., 
\[
F: \Hom(X,Y) \to \Hom(F(X), F(Y))
\]
is injective for any objects $X$, $Y$ in $\cC$. 
\end{Lemma}

\begin{proof}
Since $F$ is a C*-functor, it suffices to consider the case $X = Y$. 
From the above definition of $F(f)$, 
given $s \in S$ and $0 \not= x \in {}_1R_s$, 
the subspace 
$\displaystyle 
\begin{bmatrix}
1X\\ s
\end{bmatrix}
\otimes x\o^{1/2}$ 
of $F(X)$ is invariant under the action of $F(\End(X))$ and is equivalent 
to the obvious irreducible representation of $\End(X)$ on the vector space 
$\displaystyle
\begin{bmatrix}
X\\ s
\end{bmatrix}$. Since the family $\{ \Hom(s,X) \}_{s \in \supp(X)}$ 
is a complete set of irreducible representations of $\End(X)$, 
we conclude that $f \mapsto F(f)$ is a faithful representation of 
$\End(X)$. 
\end{proof}

\section{Free Product Functors}
\bigskip
%\noindent\underline{Operator-Valued Inner Products}
In this section, we shall reveal the monoidal structure 
of free product functors towards our realization theorem. 

The following formula on operator-valued inner products is 
frequently used to identify relative tensor products. 

\begin{Lemma}
Each element in ${}_AX_A$ is right $N$- and left $N$-bounded at the same 
time and the associated operator valued-inner products (relative to 
the normalized trace of $N$) is given by 
\[
[\x,\y]_N = [\x,\y]_A 
\qquad\text{and}
\qquad
{}_N[\x,\y] = {}_A[\x,\y].
\]
Theorefore, for $a$, $b \in N$, we have 
\[
[\x a,\y b]_N = a^*[\x,\y]_Ab, 
\qquad 
{}_N[a\x, b\y] = a{}_A[\x,\y]b^*. 
\]
\end{Lemma}

\begin{proof}
Let $\cN \subset N$ be the algebraic free product subalgebra. 
For $x \in \cN$ of the form 
\[
x = aqb_1q_1b_2q_2 \dots, 
\]
with $a \in A$, $q \in Q$, $q_i \in Q^\circ$ and $b_i \in B^\circ$ 
(here $Q^\circ = \ker \t_Q$ and $B^\circ = \ker E^B_A$), we have 
\[
\x x = \x a \otimes q\t_Q^{1/2} \otimes_A b_1\t^{1/2} \otimes_A 
\dots
\]
($\otimes_A$ is with respect to the canonical trace $\t$ on $A$) and hence 
\begin{align*}
\| \x x\|^2 &= 
(\t^{1/2} \otimes q\t_Q^{1/2} \otimes_A b_1\t^{1/2} \otimes_A 
q_1\t_Q^{1/2} \otimes \dots|\\
&\qquad
[\x a,\x a]_A \t^{1/2} \otimes_A q\t_Q^{1/2} \otimes b_1\t^{1/2} 
\otimes_A q_1\t_Q^{1/2} \otimes \dots)\\
&= (\t_N^{1/2}x| [\x,\x]_A \t_N^{1/2}x)
\end{align*}
shows that $[\x,\x]_N = [\x,\x]_A \in A \subset \cN$. 
\end{proof}

We now show the multiplicativity of free product functors. 

\begin{Lemma}
Let $X$ and $Y$ be objects in $\cC$. Then we can define 
an $N$-$N$ linear unitary map 
$m_{X,Y}: F(X)\otimes_NF(Y) \to F(X\otimes Y)$ 
by the formula 
\[
m_{X,Y}\bigl( 
(\x\otimes x\o^{1/2})\otimes_N(\y\otimes y\o^{1/2})
\bigr) 
= (\x\otimes 1_Y)\y \otimes xy\o^{1/2}, 
\]
where 
$\displaystyle
\x \in 
\begin{bmatrix}
sX\\ s'
\end{bmatrix}$, 
$\displaystyle 
\y \in 
\begin{bmatrix}
tY\\ t'
\end{bmatrix}$, 
$x \in {}_sR_{s'}$ and $y \in {}_tR_{t'}$ 
with $s$, $s'$, $t$ and $t' \in S$. 
\end{Lemma}

\begin{proof}
Set 
$\displaystyle 
\z = (\x\otimes 1_Y)\y \in 
\begin{bmatrix}
sXY\\ t'
\end{bmatrix}$. 
Then, from the operator-valued inner product formula,
\begin{align*}
\| (\x\otimes x\o^{1/2})\otimes_N (\y\otimes y\o^{1/2}) \|^2 
&= (\y\otimes y\o^{1/2}| 
[\x\otimes x\o^{1/2}, \x\otimes x\o^{1/2}]_A (\y\otimes y\o^{1/2}))\\ 
&= (\y\otimes y\o^{1/2}|(\x^*\x\otimes x^*x) (\y\otimes y\o^{1/2}))\\
&= (\z\otimes xy\o^{1/2}|\z\otimes xy\o^{1/2}).
\end{align*}
Thus we can define an isometry of ${}_AX\otimes_N Y_A$ to 
$F(XY)$ by the same formula as $m_{X,Y}$, 
which is obviously $A$-$A$ linear and has the range ${}_AXY_A$. 

Now, for $a$, $b \in \cN$, writing 
$a(\x\otimes x\o^{1/2}) = \a\otimes_A(\x\otimes x\o^{1/2})$ 
and 
$(\y\otimes y\o^{1/2})b = (\y\otimes y\o^{1/2})\otimes_A \b$ 
with $\a$, $\b$ elements in alternate $A$-tensor products of 
$L^2(Q)^\circ$ and $L^2(B)^\circ$, we have 
\begin{align*}
&\| a(\x\otimes x\o^{1/2})\otimes_N (\y\otimes y\o^{1/2})n \|^2\\
&\qquad= \| \a\otimes_A(\x\otimes x\o^{1/2})\otimes_N 
(\y\otimes y\o^{1/2})\otimes_A\b \|^2\\
&\qquad= ((\x\otimes x\o^{1/2})\otimes_N (\y\otimes y\o^{1/2}) | 
[\a,\a]_A (\x\otimes x\o^{1/2})\otimes_N (\y\otimes y\o^{1/2}) 
{}_A[\b,\b])\\
&\qquad= (\z\otimes xy\o^{1/2}| [\a,\a]_A (\z\otimes xy\o^{1/2}) 
{}_A[\b,\b])\\
\intertext{(here we apply the isometry just defined)}
&\qquad= \| \a\otimes_A(\z\otimes xy\o^{1/2})\otimes_A\b \|^2\\
&\qquad= \| a(\z\otimes xy\o^{1/2})b \|^2. 
\end{align*}
\end{proof}

The unitary map $m_{X,Y}$ is apparently natural in variables $X$, $Y$. 

\begin{Lemma}
The natural family $\{ m_{X,Y} \}$ is associative: 
\[
m_{XY,Z}(m_{X,Y}\otimes 1_{F(Z)}) = 
m_{X,YZ}(1_{F(X)}\otimes m_{Y,Z}) 
\]
for objects $X$, $Y$ and $Z$ in $\cC$. 
\end{Lemma}

\begin{proof}
We first show the associativity on the subspace 
${}_AX_A \otimes_N {}_AY_A \otimes_N {}_AZ_A$, which is 
checked by 
\begin{align*}
&m_{XY,Z}( m_{X,Y}((\x\otimes x\o^{1/2})\otimes_N(\y\otimes y\o^{1/2})) 
\otimes_N (\z\otimes z\o^{1/2}) )\\
&\qquad= m_{XY,Z}( ((\x\otimes 1_Y)\y \otimes xy\o^{1/2}) \otimes_N 
(\z\otimes z\o^{1/2}) )\\
&\qquad= (\x\otimes 1_{YZ})(\y\otimes 1_Z)\z \otimes xyz\o^{1/2}\\
&\qquad= m_{X,YZ}( 
(\x\otimes x\o^{1/2})\otimes_N 
m_{Y,Z}( (\y\otimes y\o^{1/2}) \otimes_N (\z\otimes z\o^{1/2}) ) ).
\end{align*}

Now the associativity is extended to the whole space as follows: 
By the density of $N{}_AX_A$ and ${}_AX_AN$ in $F(X)$, we see that 
\begin{align*}
\overline{ 
N({}_AX_A \otimes_N {}_AY_A \otimes_N {}_AZ_A)N }
&= \overline{ F(X) \otimes_N {}_AY_A \otimes_N {}_AZ_AN }\\
&= \overline{ F(X) \otimes_N \cN {}_AY_A \otimes_N 
{}_AZ_AN }\\
&= F(X)\otimes_N F(Y)\otimes_N F(Z). 
\end{align*}
Since $m_{X,Y}$'s are $N$-$N$ linear, the above density 
ensures the overall validity of the desired associativity. 
\end{proof}

So far, we have checked that the free product functor $F$ is monoidal 
with multiplicativity given by the family $\{ m_{X,Y} \}$. 

We next show that the monoidal functor $F$ preserves conjugations. 

\begin{Lemma}
The natural isomorphism ${}_AX^*_A \to ({}_AX_A)^*$ is extended to 
a unitary map $s_X: F(X^*) \to F(X)^*$ so that it intertwines 
$N$-$N$ actions. 
\end{Lemma}

\begin{proof}
Recall that the $A$-$A$ linear unitary map 
${}_AX^*_A \to ({}_AX_A)^*$ is defined by 
\[
\x^\star \otimes x^* \o^{1/2} 
\mapsto (\x\otimes x\o^{1/2})^*,
\qquad 
\x \in 
\begin{bmatrix}
sX\\ t
\end{bmatrix}, \ 
x \in {}_sR_t, 
\]
where 
$\displaystyle 
\x^\star \in 
\begin{bmatrix}
tX^*\\ s
\end{bmatrix}$
denotes the Frobenius transform of 
$\displaystyle 
\x^* \in 
\begin{bmatrix}
t\\ sX
\end{bmatrix}$. 

Let $a$, $b \in \cN$ and write 
$a(\x\otimes x\o^{1/2})^*b 
= \a\otimes_A (\x\otimes x\o^{1/2})^* \otimes_A\b$ as before. 
Then we have 
\begin{align*}
\| a(\x\otimes x\o^{1/2})^*b \|^2 
&= \| \a\otimes_A(\x\otimes x\o^{1/2})^*\otimes_A \b \|^2\\
&= ( (\x\otimes x\o^{1/2})^* | 
[\a,\a]_A (\x\otimes x\o^{1/2})^* {}_A[\b,\b] )\\
&= (\x^\star \otimes x^*\o^{1/2} | 
[\a,\a]_A (\x^\star \otimes x^*\o^{1/2}) {}_A[\b,\b])\\
&= \| \a\otimes_A (\x^\star \otimes x^*\o^{1/2})\otimes_A \b \|^2\\
&= \| a (\x^\star \otimes x^*\p^{1/2}) b \|^2.
\end{align*}
Thus 
\[
a(\x^\star \otimes x^*\o^{1/2})b \mapsto 
a(\x \otimes x\o^{1/2})^*b 
\]
defines a unitary map $s_X: F(X^*) \to F(X)^*$. 
\end{proof}

%\begin{figure}[htbp]
%%\vspace*{-0.5cm}
%%\hspace*{3cm}
%\includegraphics{fig4.eps}
%%\vspace*{3cm}
%\caption[]{}
%\end{figure}

\begin{Lemma}
The family $\{ s_X \}$ is natural in $X$: 
The diagram
\[
\begin{CD}
F(X^*) @>{s_X}>> F(X)^*\\
@A{F({}^tf)}AA @AA{{}^tF(f)}A\\
F(Y^*) @>>{s_Y}> F(Y)^*
\end{CD}
\]
commutes for $f: X \to Y$. 
\end{Lemma}

\begin{proof}
By the naturality of Frobenius duality, we have 
$(1\otimes {}^tf)\y^\star = ((1\otimes f^*)\y)^\star$, 
which is utilized in the following way: 
\begin{align*}
s_X F({}^tf) (\y^\star \otimes y^*\o^{1/2}) 
&= s_X((1\otimes {}^tf)\y^\star \otimes y^*\o^{1/2})\\ 
&= s_X(((1\otimes f^*)\y)^\star \otimes y^*\o^{1/2})\\ 
&= ( (1\otimes f^*)\y \otimes y\o^{1/2} )^*\\
&= \bigl( F(f^*)(\y\otimes y\o^{1/2}) \bigr)^*\\ 
&= \bigl( F(f)^*(\y\otimes y\o^{1/2}) \bigr)^*\\ 
&= {}^tF(f) (\y\otimes y\o^{1/2})^*.
\end{align*}
\end{proof}

%\begin{figure}[htbp]
%%\vspace*{-0.5cm}
%%\hspace*{-9cm}
%\includegraphics{fig5.eps}
%%\vspace*{4cm}
%\caption[]{}
%\end{figure}

\begin{Lemma}
The family $\{ s_X \}$ is multiplicative: 
\[
\begin{CD}
F(Y^*)\otimes_N F(X^*) @>{s_Y\otimes s_X}>> F(Y)^*\otimes_NF(X)^*\\
@V{m_{Y^*,X^*}}VV @AA{{}^tm_{X,Y}}A\\
F(Y^*\otimes X^*) @>>{s_{XY}}> F(X\otimes Y)^*
\end{CD}\quad.
\]
\end{Lemma}

\begin{proof}
By ${}^tm_{X,Y}^{-1} = \overline{m_{X,Y}}$, 
we shall verify the equivalent relation 
$\overline{m_{X,Y}} (s_Y\otimes s_X) = s_{XY} m_{Y^*,X^*}$, 
which is checked on the $N$-cyclic subspace in the following way:
\begin{align*}
\overline{m_{X,Y}} (s_Y\otimes s_X) 
&\bigl( (\y^\star\otimes y^*\o^{1/2}) \otimes_N (\x^\star\otimes
x^*\o^{1/2}) \bigr)\\
&= \overline{m_{X,Y}} ( (\y\otimes y\o^{1/2})^* \otimes_N 
(\x\otimes x\o^{1/2})^* )\\
&= \overline{m_{X,Y}} ( (\x\otimes x\o^{1/2}) \otimes_N 
(\y\otimes y\o^{1/2}) )^*\\
&= \bigl( 
m_{X,Y}((\x\otimes x\o^{1/2})\otimes_N(\y\otimes y\o^{1/2})) 
\bigr)^*\\
&= \bigl( (\x\otimes 1_Y)\y\otimes xy\o^{1/2} \bigr)^*,
\end{align*}
which coincides with 
\[
s_{XY} m_{Y^*,X^*} ( (\y^\star\otimes y^*\o^{1/2}) 
\otimes_N (\x^\star\otimes x^*\o^{1/2}) ) 
\]
if we use
\[
(\y^\star\otimes 1_{X^*})\x^\star = 
\bigl( (\x\otimes 1_Y)\y \bigr)^\star, 
\]
which is an easy consequence of Frobenius transforms. 
\end{proof}

%\begin{figure}[htbp]
%\vspace*{1cm}
%%\hspace*{-11cm}
%\includegraphics{fig6.eps}
%%\vspace*{5cm}
%\caption[]{}
%\end{figure}

\begin{Lemma}
The family $\{ s_X \}$ is compatible with duality: 
Two unitary maps 
$s_{X^*}: F(X^{**}) \to F(X^*)^*$
and
${}^ts_X: F(X)^{**} \to F(X^*)^*$ are the same 
if we apply the identification
$F(X^{**}) = F(X) = F(X)^{**}$. 
\end{Lemma}

\begin{proof}
We compute as 
\begin{align*}
\overline{s_X}s_{X^*} (\x\otimes x\o^{1/2}) 
&= \overline{s_X} (\x^\star\otimes x^*\o^{1/2})^* 
= \bigl( s_X(\x^\star\otimes x^*\o^{1/2}) \bigr)^*\\
&= (\x\otimes x\o^{1/2})^{**}
= \x\otimes x\o^{1/2},
\end{align*}
which implies $\overline{s_X} s_{X^*} = 1_{F(X)}$ by the density
argument and we are done as ${}^ts_X^{-1} = \overline{s_X}$. 
\end{proof}

We can now conclude that the free product functor is *-monoidal 
with the accompanied isomorphisms given by $\{ s_X \}$. 

From the properties of the free product functor established so far, 
we find that, for each object $K$ in $\cC$, 
$F(K)^*$ is a dual object of $F(K)$ in 
the C*-tensor category of $N$-$N$ bimodules. 
Since $N$ is a factor, the unit object $L^2(N)$ is irreducible 
and hence the Longo-Roberts finite-dimensionality theorem 
can be applied to show that $\End(F(K))$ is finite-dimensional 
and at the same time each $F(K)$ has finite Jones index by 
the rigidity characterization of index-finiteness (Proposition~1.9). 
Moreover we know the following weaker version of Frobenius 
reciprocity 
\[
\End(F(K)) \cong \Hom(L^2(N), F(K)\otimes_N F(K)^*).
\]
By the natural identification 
\[
F(K)\otimes_N F(K)^* \cong F(KK^*) \cong 
L^2\bigl( (Q\otimes A(K))*_{A(K)}B(K) \bigr), 
\]
the above fact implies that $\End(F(K))$ is isomorphic to 
the set of $N$-central vectors in 
$L^2((Q\otimes A(K))*_{A(K)}B(K))$, which is further isomorphic to 
the relative commutant 
\[
N' \cap \bigl( (Q\otimes A(K))*_{A(K)}B(K) \bigr). 
\]
Now the lemma below shows that this is in fact isomorphic to 
$\End(K)$, whence 
\[
\End (K) \ni f \mapsto F(f) \in \End(F(K))
\]
is a surjective isomorphism for any $K$ and then reducing it to
off-diagonal corners, the functor $F$ is found to be fully 
faithful: 
\[
F: \Hom(X,Y) \to \Hom(F(X),F(Y))
\]
is a surjective isomorphism for any objects $X$, $Y$ in $\cC$. 

\begin{Lemma}
We have 
\[
N' \cap \bigl( (Q\otimes A(K))*_{A(K)}B(K) \bigr) = \End(K). 
\]
\end{Lemma}

\begin{proof}
By the ergodicity property of free products 
(Theorem~4.1 in \cite{P1}), we have 
\[
Q' \cap 
\bigl( (Q\otimes A(K))*_{A(K)}B(K) \bigr) = Q \otimes A(K) 
\]
and hence 
\[
N' \cap 
\bigl( (Q\otimes A(K))*_{A(K)}B(K) \bigr) = B' \cap A(K). 
\]
%$Q^\circ\otimes A(K)$ and $B^\circ$ are free.

Taking into account the commutativity with elements in $B$ of the form 
$1_s \otimes x$ with $s \in S^2$ and $x \in {}_sR_s$, 
we first reduce an element in $B'\cap A(K)$ to the form 
\[
\s = \bigoplus_{s \in S} \s_s \otimes e_s 
\qquad\text{with}\ 
\s_s \in 
\begin{bmatrix}
sK\\ sK
\end{bmatrix}. 
\]

Now the commutativity with the Jones projection 
$\displaystyle d(s)^{-1} \e_{s^*}^*\e_{s^*} \in 
\begin{bmatrix}
s^*s\\ s^*s
\end{bmatrix}$ 
reveals that $\s_s$ is further restricted to 
\[
\s_s = 1_s\otimes f_s 
\qquad\text{with}\quad 
f_s \in \End(K). 
\]

To see the independance of $f_s$ on $s \in S$, 
let $s$, $t \in S$ and 
$\displaystyle \r = \e_{s^*}^*\e_{t^*} \in 
\begin{bmatrix}
s^*s\\ t^*t
\end{bmatrix}$. 
Then the commutativity of $\s$ with 
$\r\otimes {}_{s^*s}R_{t^*t} \subset B$ 
means 
\[
(\r\otimes f_s) \otimes x = (\r\otimes f_t) \otimes x
\]
for $x \in {}_{s^*s}R_{t^*t}$. Since ${}_{s^*s}R_{t^*t} \not= \{0\}$ 
and $\rho \not= 0$, this implies $f_s = f_t$ for any $s$, $t \in S$.

Thus, letting $f = f_s \in \End(K)$, we see that $\s = f$ belongs to 
$\End(K)$. 
\end{proof}

\begin{Theorem}
Let $\cC$ be a C*-tensor category with Frobenius duality. 
Then the free product functor $F$ gives a fully faithful realization 
of $\cC$ as that of $N$-$N$ bimodules of finite index. 
Moreover, $F$ preserves Frobenius dualities as well as conjugations. 
\end{Theorem}

\begin{proof}
Non-trivial is the equality $F(\e_X) = \e_{F(X)}$. 
By the additivity of $\{ \e_X \}$ and $\{ \e_{F(X)} \}$, 
we may restrict ourselves to 
the case of a simple $X$, which is checked by showing 
the positivity of $F(\e_X)$: 
\begin{align*}
F(\e_X) 
\bigl( 
(\x\otimes x\o^{1/2})\otimes_N (\x\otimes x\o^{1/2})^* 
\bigr)
&= F(\e_X) 
\bigl( 
(\x\otimes 1_{X^*})\x^\star \otimes xx^*\o^{1/2} 
\bigr)\\
&= (1\otimes \e_X)(\x\otimes 1_{X^*})\x^\star \otimes xx^*\o^{1/2}\\
&= \langle \x\x^* \rangle_X \otimes xx^*\o^{1/2},
\end{align*}
which obviously belongs to the positive cone of $L^2(N)$. 
\end{proof}

%\begin{figure}[htbp]
%%\vspace*{-0.5cm}
%%\hspace*{-4.5cm}
%\includegraphics{fig7.eps}
%%\vspace*{4.5cm}
%\caption[]{}
%\end{figure}

\section{Comments}

If we apply our main theorem to the Tannaka dual of a compact quantum 
group $G$ (see \cite{Y2} for the discussion of Frobenius duality), 
the resulting bimodule realization can be interpreted as 
giving a kind of coaction of $G$ on the finite 
factor $N$ (the reference von Neumann algebra of bimodules), 
which is minimal because the realization is fully faithful (\cite{Y3}). 
The crossed product algebra $M$ is then a factor and 
the ``dual'' action of $G$ on $M$ is minimal in the sense that 
$M'\cap (M\rtimes G) = \C 1$ by Takesaki's duality 
(see \cite[\S 8]{HY} for details). 

Thus we can recover the result due to Y.~Ueda \cite{U1}. 
Note that his construction is more direct than ours:
the minimal action is described as the free product of a faithful action 
and a trivial action, whereas subfactors corresponding to our 
bimodule realization are constructed according to 
the method of A.~Wassermann (\cite{U2}). 
In this sense, our construction is a kind of reverse process of Ueda's
although it is not clear at present whether they really give 
the same subfactors.

%The computation of amalgamated free product factors due to 
%Radulescu and Dykema strongly suggest that our factors $N(Q)$ will be all 
%isomorphic as long as $Q$ is continuous. 

If we apply our realization to the bicategory generated by 
a single object $X$, then we obtain a realization of 
the associated standard lattice as higher relative commutants 
of the subfactor $N \subset \End(X_N)$, recovering the result 
due to S.~Popa. It is therefore interesting to make clear 
the relationship between the Popa's realization and ours. 
Roughly speaking, our free product factors are amplified in its
ingredients by 
the AFD $\text{II}_1$-factor 
compared with the Popa's subfactors but no clues to explicit connections 
presently. 

Finally we remark here that our realization theorem cannot be obtained by just 
reformulating the Popa's result: 
Firstly, although it is straightforward to produce standard lattices 
from tensor categories (with Frobenius duality), the reverse 
implication is highly non-trivial and is not established yet
(cf.~\cite{J2}). Secondly, standard lattices should correspond to 
singly generated tensor categories and there would be no good way to 
approximate tensor categories by their subcategories especially 
in the realization problem (a kind of cohomological adjustments 
are needed among approximating realizations, which seem hopeless).

\end{document}